\documentclass{amsart}
\def\B{{\mathcal B}}
\def\H{{\mathcal H}}
\def\Le{{\mathcal L}}

\def\Pe{{\mathcal P}}
\def\X{{\mathcal X}}
\def\Y{{\mathcal Y}}
\def\V{{\kern-.5pt\mathcal V}}
\def\W{{\mathcal W}}
\def\Z{{\mathcal Z}}

\def\CC{{\mathbb C\kern.5pt}}
\def\FF{{\mathbb F\kern.5pt}}
\def\LL{{\kern.5pt\mathbb L\kern.5pt}}
\def\RR{{\mathbb R\kern.5pt}}

\newsymbol\Bbbk 207C
\newsymbol\diagdown 231F
\newsymbol\digamma 207A

\let\vphi=\varphi
\let\sse=\subseteq
\let\what=\widehat

\def\newmatrix#1{\null\,\vcenter{
		\baselineskip=8pt\mathsurround=-0pt\ialign{
		\hfil ${##}$
		\hfil &&
		\hfil ${##}$
		\hfil \crcr
		\mathstrut \crcr
		\noalign{\kern-\baselineskip}#1 \crcr
		\mathstrut \crcr
		\noalign{\kern-\baselineskip} \crcr }}\!}

\def\x{{\times}}
\def\noi{\noindent}
\def\ve{{_\vee\kern-1pt}}
\def\we{{_\wedge\kern-1pt}}
\def\hotimes{{\kern1pt\what\otimes}}
\def\otimesv{{\kern1pt\otimes_\ve}}
\def\otimesw{{\kern1pt\otimes_{\kern-1pt\we}}}
\def\hotimesv{{\kern1pt\what\otimes_\ve}}
\def\hotimesw{{\kern1pt\what\otimes_{\kern-1pt_\we}}}
\def\smallfrac#1#2{{\textstyle{\frac{#1}{#2}}}}

\def\emap{\hbox to25pt{\rightarrowfill}}

\def\nmap{\Big\uparrow}
\def\diagdownBOX{\hbox{$\diagdown$}}
\def\searrowBOX{\hbox{\hglue6.5pt$\searrow$}}
\def\semap{\vbox{\offinterlineskip
\diagdownBOX\vglue-1pt\searrowBOX\vglue-6pt}}

\def\theorem #1{\vskip6pt\noi{\bf{Theorem} #1$.$}}

\def\corollary #1{\vskip6pt\noi{\bf{Corollary} #1$.$}}
\def\proposition #1{\vskip6pt\noi{\bf{Proposition} #1$.$}}
\def\claim #1{\vskip6pt\noi{\it{\kern-.5pt Claim} #1$.$}}
\def\remark #1{\vskip6pt\noi{\bf{Remark} #1$.$}}

\begin{document}

\vglue-55pt\noi
\hfill{\it New York Journal of Mathematics}\/,
{\bf 28} (2022) 1656--1666

\vglue30pt
\title
{Restriction of Uniform Crossnorms}
\author{Carlos S. Kubrusly}
\address{Catholic University of Rio de Janeiro, Brazil}
\email{carlos@ele.puc-rio.br}
\date{September 20, 2022 (minor corrections: February 3, 2023)}

\begin{abstract}
Consider the space $\B[\X\otimes_\alpha\!\Y]$ of operators on the tensors
product ${\X\otimes_\alpha\!\Y}$ of normed spaces $\X$ and $\Y$ equipped with
a uniform crossnorm ${\|\cdot\|_\alpha}$$.$ Take the induced uniform norm
on $\B[\X\otimes_\alpha\!\Y]$ and consider its restriction to the tensor
product ${\B[\X]\otimes\B[\Y]}$ of the algebra of operators $\B[\X]$ and
$\B[\Y].$ It is proved that such a restriction is a reasonable crossnorm
on ${\B[\X]\otimes\B[\Y]}$.
\end{abstract}
\vskip-10pt

\maketitle

\section{Introduction}

The paper deals with a special tensor norm on the tensor product of a pair of
spaces of bounded linear transformations; in particular, of a pair of spaces
of operators$.$ We avoid the term ``operator space'' in this note, since the
term has already been consecrated to define a theory of certain subspaces of
the algebra of Hilbert-space operators $\B[\H]$, which can be thought of as
object of a category in the realm of $C^*$-algebras (see, e.g., \cite{BP},
\cite[Definition 1.2]{Pis1}, \cite{Arv}, \cite[Definition 1.1]{Pis2})$.$ On
the contrary, our aim in this note is much less ambitious$.$ Let $\B[\X]$ be
the normed algebra of all operators on a normed space $\X$, consider the
tensor product ${\B[\X]\otimes\B[\Y]}.$ It is shown that this is included in
the normed algebra ${\B[\X\otimes_\alpha\!\Y]}$, were ${\X\otimes_\alpha\!\Y}$
stands for the tensor product ${\X\otimes\Y}$ of normed spaces $\X$ and $\Y$
equipped with a uniform crossnorm ${\|\cdot\|_\alpha}.$ We give an elementary
proof that the induced uniform operator norm on ${\B[\X\otimes_\alpha\!\Y]}$,
when restricted to ${\B[\X]\otimes\B[\Y]}$, acts as a reasonable crossnorm on
the tensor product ${\B[\X]\otimes\B[\Y]}$.

\vskip6pt
All terms and notation used above will be defined here in due course$.$ The
paper is organized into four more sections$.$ Basic propositions, including
notation and terminology, are summarized in Section 2$.$ Supplementary results
on uniform crossnorms required in the sequel are considered in Section 3$.$
The main theorem is proved in Section 4$.$ A discussion on equivalent uniform
crossnorms, in light of the outcome of Section 4, closes the paper in
Section 5.

\section{Auxiliary Results}

All linear spaces in this paper are over the same scalar field $\FF$, which is
either $\RR$ or $\CC.$ The algebraic {\it tensor product}\/ of linear spaces
$\X$ and $\Y$ is a linear space ${\X\otimes\Y}$ for which there is a bilinear
map ${\theta\!:\X\x\Y\to\X\otimes\Y}$ (called the {\it natural bilinear map}\/
associated with ${\X\otimes\Y}$) whose range spans ${\X\otimes\Y}$ with the
following additional (universal) property: for every bilinear map
${\phi\!:\X\x\Y\to\Z}$ into any linear space $\Z$ there exists a (unique)
linear transformation ${\Phi\!:\X\otimes\Y\to\Z}$ for which the diagram
$$
\newmatrix{
\X\x\Y & \kern2pt\buildrel\phi\over\emap & \kern-1pt\Z                   \cr
       &                                 &                               \cr
       & \kern-3pt_\theta\kern-3pt\semap & \kern4pt\nmap\scriptstyle\Phi \cr
       &                                 & \phantom{;}                   \cr
       &                                 & \kern-2pt\X\otimes\Y          \cr}
$$
commutes$.$ A tensor product space ${\X\otimes\Y}$ exists for every pair of
linear spaces ${(\X,\Y)}$ and is unique up to isomorphisms$.$ Set
${x\otimes y=\theta(x,y)\in\X\otimes\Y}$ for each ${(x,y)\in\X\x\Y}.$ These
are the {\it single tensors}\/$.$ An arbitrary element $\digamma$ in the
linear space ${\X\otimes\Y}$ is a finite sum ${\sum_ix_i\otimes y_i}$ of
single tensors, and the representation of $\digamma\!={\sum_ix_i\otimes y_i}$
as a finite sum of single tensors is not unique$.$ (For an exposition on
algebraic tensor products see, e.g., \cite{Kub}$.)$ If $\X$ and $\Y$ are
linear spaces, then ${\Le[\X,\Y]}$ denotes the linear space of all linear
transformations of $\X$ into $\Y.$ Let ${\X,\Y,\V,\W}$ be linear spaces
and consider the tensor product spaces ${\X\otimes\Y}$ and ${\V\otimes\W}.$
Take a pair of linear transformations
${A\in\Le[\X,\V]}$ and ${B\in\Le[\Y,\W]}$ and set
$$
(A\otimes B){\sum}_ix_i\otimes y_i={\sum}_iAx_i\otimes By_i
$$
in ${\V\otimes\W}$ for every ${\digamma\!=\sum_ix_i\otimes y_i}$ in
${\X\otimes\Y}.$ This defines a linear transformation
${A\otimes B\in\Le[\X\otimes\Y,\V\otimes\W]}$ of the linear space
${\X\otimes\Y}$ into the linear space ${\V\otimes\W}$, referred to as the
{\it tensor product of the linear transformations}\/ $A$ and $B$, which is
such that $(A\otimes B)(\digamma)$ does not depend on the representation
${{\sum}_ix_i\!\otimes\!y_i}$ of ${\digamma\!\in\X\otimes\Y}$ (see, e.g.,
\cite[Proposition 3.6]{Kub})$.$ Consider the linear spaces ${\Le[\X,\V]}$
and ${\Le[\Y,\W]}$ and let $\LL$ be an arbitrary element of
${\Le[\X,\V]\otimes\Le[\Y,\W]}$ so that ${\LL=\sum_jA_j\otimes B_j}$ is a
finite sum of single tensors ${A_j\otimes B_j\in\Le[\X,\V]\otimes\Le[\Y,\W]}$,
and therefore
$$
\Le[\X,\V]\otimes\Le[\Y,\W]\sse\Le[\X\otimes\Y,\V\otimes\W].
$$
\vskip-2pt

\vskip6pt
{From} now on suppose $\X$ and $\Y$ are normed spaces$.$ Let ${\B[\X,\Y]}$ be
the normed space of all bounded linear transformations of $\X$ into $\Y$
equipped with its standard induced uniform norm, and let ${\X^*\!=\B[\X,\FF]}$
denote the dual of $\X.$ Let ${x\otimes y}$ and ${f\otimes g}$ be single
tensors in the tensor product spaces ${\X\otimes\Y}$ and ${\X^*\otimes\Y^*}.$
A norm ${\|\cdot\|_\alpha}\!$ on ${\X\otimes\Y}$ is a {\it reasonable
crossnorm}\/ if for every ${x\in\X}$, ${y\in\Y}$, \hbox{${f\in\X^*}\!$,
${g\in\Y^*}\!$},
\begin{description}
\item {$\kern-12pt$\rm(a)}
$\;\|x\otimes y\|_\alpha\kern-1pt\le\|x\|\,\|y\|$,
\vskip4pt
\item {$\kern-12pt$\rm(b)}
$\;{f\otimes g}$ lies in $(\X\otimes\Y)^*\!$, $\,$and
${\|f\otimes g\|_{*\alpha}\kern-1pt\le\|f\|\,\|g\|}\;$ (where
${\|\cdot\|_{*\alpha}}$ is the norm on the dual ${(\X\otimes\Y)^*}$ when
${\X\otimes\Y}$ is equipped with the norm ${\|\cdot\|_\alpha}$), so that
\end{description}
$$
\X^*\otimes\Y^*\sse(\X\otimes\Y)^*.
$$
Actually, ${\|x\otimes y\|_\alpha\kern-1pt=\|x\|\,\|y\|}$ and
${\|f\otimes g\|_{*\alpha}\kern-1pt=\|f\|\,\|g\|}$ whenever
${\|\cdot\|_\alpha}\kern-1pt$ is a reasonable crossnorm (see, e.g.,
\cite[Proposition 1.1.1]{DFS})$.$ Two especial reasonable crossnorms on
${\X\otimes\Y}$ are the {\it injective}\/ ${\|\cdot\|_\ve}$ and
{\it projective}\/ ${\|\cdot\|_\we}$ norms,
$$
\|\digamma\|_{_\vee\kern-1pt}
=\sup_{\|f\|\le1,\,\|g\|\le1,\;f\in\X^*\!,\,g\in\Y^*}
\Big|{\sum}_if(x_i)\,g(y_i)\Big|,
$$
$$
\|\digamma\|_{_\wedge\kern-1pt}
=\inf_{\{x_i\}_i,\,\{y_i\}_i,\;\digamma=\sum_ix_i\otimes y_i}
{\sum}_i\|x_i\|\,\|y_i\|,
$$
for every ${\digamma\!=\kern-1pt\sum_i\kern-1ptx_i\otimes y_i}$ (the infimum
is taken over all representations of
${\digamma\!\in\kern-1pt\X\otimes\kern-1pt\Y}).$

\vskip3pt
\proposition{2.1}
{\it A norm\/ ${\|\cdot\|_\alpha}\kern-1pt$ on\/
${\X\otimes\kern-1pt\Y}\kern-1pt$ is a reasonable crossnorm
\hbox{if and only if}}
$$
\|\digamma\|_\ve
\le\|\digamma\|_\alpha
\le\|\digamma\|_\we
\quad\;\hbox{\it for every}\;\quad
\digamma\in\X\otimes\Y.
$$

\proof
See, e.g., \cite[Proposition 6.1]{Rya}. \qed

\vskip9pt
Let ${\X\otimes_\alpha\!\Y=(\X\otimes\Y,\|\cdot\|_\alpha)}$ denote the tensor
product space of normed spaces equipped with a norm ${\|\cdot\|_\alpha}$ on
${\X\otimes\Y}$, which is not necessarily complete even if
${\|\cdot\|_\alpha}$ is a reasonable crossnorm and $\X$ and $\Y$ are Banach
spaces$.$ Their completion is denoted by ${\X\hotimes_\alpha\Y}$ (same
notation ${\|\cdot\|_\alpha}$ for the extended norm on
${\X\hotimes_\alpha\Y}$)$.$ In particular, ${\X\hotimesv\Y}$ and
${\X\hotimesw\Y}$ are referred to as the {\it injective and projective tensor
products}\/$.$ For the theory of Banach space ${\X\hotimes_\alpha\Y}$
(including ${\X\hotimesv\Y}$ and ${\X\hotimesw\Y}$) see, e.g.,
\cite[Chapters 15 and 16]{Jar}, \cite[Section 12]{DF},
\cite[Section 6.1]{Rya}, \hbox{\cite[Sections 1.1 and 1.2]{DFS}.}

\vskip6pt
Recall that ${\X^*\otimes\Y^*\sse(\X\otimes_\alpha\!\Y)^*}.$ When restricted
to ${\X^*\otimes\Y^*}$ the norm ${\|\cdot\|_{*\alpha}}$ on
${(\X\otimes_\alpha\!\Y)^*}$ is a reasonable crossnorm (with respect to
${(\X^*\otimes\Y^*)^*})$.

\vskip3pt
\proposition{2.2}
{\it Let\/ ${\|\cdot\|_\alpha}\!$ be a reasonable crossnorm on a tensor
product space\/ ${\X\otimes\Y}$ of normed spaces\/ $\X$ and\/ $\Y$, and take
the dual\/ ${(\X\otimes_\alpha\!\Y)^*}$ of\/
${\X\otimes_\alpha\!\Y}={(\X\otimes\Y,\|\cdot\|_\alpha)}.$ $\!$When restricted
to\/ ${\X^*\otimes\Y^*}\!$ the norm\/ ${\|\cdot\|_{*\alpha}}$ on\/
${(\X\otimes_\alpha\!\Y)^*}$ is a reasonable crossnorm on}\/
${\X^*\otimes\Y^*}$.

\proof
See, e.g., \cite[Proposition 1.1.2]{DFS}. \qed

\vskip9pt
The purpose of this paper is to extend the above (nontrivial) result to the
case where ${\X^*\!=\B[\X,\FF]}$, $\,{\Y^*\!=\B[\Y,\FF]}$,
$\,{\X^*\!\otimes\Y^*\!=\B[\X,\FF]\otimes\B[\Y,\FF]}$ and
${(\X\otimes_\alpha\!\Y)^*}\!={\B[\X\otimes_\alpha\!\Y,\FF]}$ are replaced by
(extended to) $\B[\X,\V]$, $\,\B[\Y,\W]$, ${\B[\X,\V]\otimes\B[\Y,\W]}$, and
${\B[\X\otimes_\alpha\!\Y,\V\otimes_\alpha\!\W]}$, respectively, for arbitrary
normed spaces ${\X,\Y,\V,\W}.$ This will be done in Section 4 (Theorem 4.1).

\section{Uniform Crossnorms}

If $\X$, $\Y$, and $\Z$ are normed spaces and if ${T\in\B[\X,\Y]}$ and
${S\in\B[\Y,\Z]},$ then ${S\kern1ptT\in\B[\X,\Z]}$ and
${\|S\kern1ptT\|\le\|S\|\,\|T\|}.$ This is a crucial property shared by the
induced uniform norm of bounded linear transformations, referred to as the
operator norm property (see, e.g.,\cite[Proposition 4.16]{EOT})$.$ Its
counterpart for the case of tensor products (rather than ordinary products)
yields the notion of uniform crossnorm.

\vskip6pt
A {\it uniform crossnorm}\/ ${\|\cdot\|_\alpha}$ is a reasonable crossnorm on
every tensor product space (of arbitrary normed spaces ${\X,\Y,\V,\W}$) such
that for every bounded linear transformations ${A\in\B[\X,\V]}$ and
${B\in\B[\Y,\W]}$ the linear tensor product transforma\-tion
${A\otimes B\!:\!\X\kern-1pt\otimes_\alpha\!\Y\to\!\V\otimes_\alpha\!\!\W}$
is bounded (i.e.,
${A\otimes B\in\B[\X\otimes_\alpha\!\Y,\V\otimes_\alpha\!\W]}\kern.5pt$) and
$$
\|A\otimes B\|\le\|A\|\,\|B\|,
\quad\hbox{equivalently,}\quad
\|A\otimes B\|=\|A\|\,\|B\|.
$$
In fact, the above equivalence holds since
\begin{eqnarray*}
\|A\|\kern1pt\|B\|
&\kern-6pt=\kern-6pt&
\sup_{\|x\|\le1}\|Ax\|\sup_{\|y\|\le1}\|By\|
\le\sup_{\|x\|\|y\|\le1}\|Ax\|\kern1pt\|By\|
=\sup_{\|x\otimes y\|_\alpha\le1}\|Ax\otimes By\|_\alpha                  \\
&\kern-6pt=\kern-6pt&
\sup_{\|x\otimes y\|_\alpha\le1}\|(A\otimes B)(x\otimes y)\|_\alpha
\le\sup_{\|\digamma\|_\alpha\le1}\|(A\otimes B)\digamma\|_\alpha
=\|A\otimes B\|,
\end{eqnarray*}
with ${\|A\!\otimes\!B\|}$, $\|A\|$, and $\|B\|$ standing for the induced
uniform norms of ${A\otimes B}$ in
${\B[\X\otimes_\alpha\!\Y,\,\V\otimes_\alpha\!\W]}$, $\,A$ in ${\B[\X,\V]}$,
and $B$ in ${\B[\Y,\W]}.$

\vskip6pt
The projective ${\|\cdot\|_\we}$ and injective
${\|\cdot\|_\ve}$ norms are uniform crossnorms (see, e.g.,
\hbox{\cite[Propositions 2.3 and 3.2]{Rya}\kern.5pt)}$.$

\vskip6pt
Let ${\|\cdot\|_{[\alpha,\alpha]}}$ denote the induced uniform norm on the
normed space of bounded linear transformations
${\B[\X\otimes_\alpha\!\Y_,\,\V\otimes_\alpha\!\W]}$ (i.e.,
${\|\cdot\|_{[\alpha,\alpha]}}$ stands for the induced uniform norm
${\|\cdot\|}$ on ${\B[\X\otimes\Y,\V\otimes\W]}$ when
${\X\kern-1pt\otimes\kern-1pt\V}$ and ${\Y\otimes\kern-1pt\W}$ are equipped
with a uniform crossnorm ${\|\cdot\|_\alpha}).$ The notation
${\|\cdot\|_{[\alpha,\alpha]}}$ highlights the dependence of ${\|\cdot\|}$ on
the uniform crossnorm ${\|\cdot\|_\alpha}$ which equips both ${\X\otimes\Y}$
and ${\V\otimes\W}.$ For instance,
$$
\|A\otimes B\|_{[\alpha,\alpha]}
=\sup_{\digamma\in\X\otimes\Y,\;\|\digamma\|_\alpha\le1}
\|(A\otimes B)\digamma\|_\alpha
=\|A\|\,\|B\|.
$$
\vskip-2pt

\vskip3pt
\proposition{3.1}
{\it If\/ ${\|\cdot\|_\alpha}$ is a uniform crossnorm, then\/
${\B[\X,\V]\otimes\B[\Y,\W]}
\sse{\B[\X\otimes_\alpha\!\Y,\,\V\otimes_\alpha\!\W]}$
so that}
$$
\B[\X,\V]\otimes_{[\alpha,\alpha]}\B[\Y,\W]
\sse\B[\X\otimes_\alpha\!\Y,\,\V\otimes_\alpha\!\W].
$$

\proof
An arbitrary element $\LL$ in the tensor product space
${\B[\X,\kern-1pt\V]\otimes\B[\Y,\kern-1pt\W]}$ is represented by a finite sum
$\LL={\sum_jA_j\otimes B_j}$ of single tensors with ${A_j\in\B[\X,\V]}$
and ${B_j\in\B[\Y,\W]}.$ If ${\|\cdot\|_\alpha}$ is a uniform crossnorm, then
each linear transformation ${A_j\otimes B_j}$ lies in
${\B[\X\otimes_\alpha\!\Y,\V\otimes_\alpha\!\W]}$ and so does
$\LL={\sum_j\!A_j\!\otimes\!B_j}.$ Therefore
$$
\B[\X,\V]\otimes\B[\Y,\W]
\sse\B[\X\otimes_\alpha\!\Y,\,\V\otimes_\alpha\!\W].
$$
Now equip ${\B[\X,\V]\otimes\B[\Y,\W]}$ with the norm
${\|\cdot\|_{[\alpha,\alpha]}}$ from
${\B[\X\otimes_\alpha\!\Y,\,\V\otimes_\alpha\!\W]}$ and set
${\B[\X,\V]\otimes_{[\alpha,\alpha]}\B[\Y,\W]
=(\B[\X,\V]\otimes\B[\Y,\W],\,\|\cdot\|_{[\alpha,\alpha]})}$.          \qed

\vskip6pt
Let ${\|\cdot\|_{*[\alpha,\alpha]}}$ be the induced uniform norm on
${\B\big[\B[\X,\V]\otimes_{[\alpha,\alpha]}\B[\Y,\W],\FF\kern.5pt\big]}$.

\vskip3pt
\proposition{3.2}
{\it If\/ ${\|\cdot\|_\alpha}$ is a uniform crossnorm, then\/
$\B[\X,\V]^*\otimes\B[\Y,\W]^*
\sse\Le\big[\B[\X,\V]\!\otimes_{[\alpha,\alpha]}\!\B[\Y,\W],\FF\kern.5pt\big].$
Moreover,
$$
\|\vphi\otimes\eta\|_{*[\alpha,\alpha]}\le\|\vphi\|\,\|\eta\|
$$
for every\/ ${\vphi\in\B[\X,\V]^*}$ and\/ ${\eta\in\B[\Y,\W]^*}.$
Consequently}\/,
$$
\B[\X,\V]^*\otimes_{*[\alpha,\alpha]}\B[\Y,\W]^*
\sse(\B[\X,\V]\otimes_{[\alpha,\alpha]}\B[\Y,\W])^*.
$$

\proof
Since ${\FF\otimes\FF\cong\FF}$ (here $\cong$ stands for algebraic
isomorphism), we get
\begin{eqnarray*}
\B[\X,\V]^*\!\otimes\B[\Y,\W]^*
\!&\kern-6pt=\kern-6pt&\!
\B\big[\B[\X,\V],\FF\kern.5pt]\kern.5pt\big]
\otimes\B\big[\B[\Y,\W],\FF\kern.5pt]\kern.5pt\big]                       \\
\!&\kern-6pt\sse\kern-6pt&\!
\Le\big[\B[\X,\V],\FF\kern.5pt]\kern.5pt\big]
\otimes\Le\big[\B[\Y,\W],\FF\kern.5pt]\kern.5pt\big]                      \\
\!&\kern-6pt\sse\kern-6pt&\!
\Le\big[\B[\X,\V]\otimes\B[\Y,\W],\FF\otimes\FF\kern.5pt\big]
\!=\Le\big[\B[\X,\V]\!\otimes_{[\alpha,\alpha]}\!\B[\Y,\W],\FF\kern.5pt\big],
\end{eqnarray*}
when ${\B[\X,\V]\otimes\B[\Y,\W]}$ is equipped with the uniform induced norm
${\|\cdot\|_{\alpha,\alpha]}}$ on
${\B[\X\otimes_\alpha\!\Y,\,\V\otimes_\alpha\!\W]}$ by Proposition 3.1$.$
Take ${\vphi\in\B[\X,\V]^*}$ and ${\eta\in\B[\Y,\W]^*}$ arbi\-tray$.$
According to the above inclusion,
${\vphi\otimes\eta
\in\Le\big[\B[\X,\V]\otimes_{[\alpha,\alpha]}\B[\Y,\W],\FF\kern.5pt\big]}.$
Set
$$
\|\vphi\otimes\eta\|_{*[\alpha,\alpha]}
=\sup_{\LL\in\B[\X,\V]\otimes\B[\Y,\W],\;\|\LL\|_{[\alpha,\alpha]}\le1}
|(\vphi\otimes\eta)\LL|.
$$
It was show in \cite[Theorem 3]{Sim} that the above supremum is not only
finite but bounded by ${\|\vphi\|\,\|\eta\|}$ for the particular case of
${\V=\X}$ and ${\W=\Y}$ when these are Banach spaces, whose extension for
arbitrary normed spaces ${\X,\Y,\V,\W}$ follows the same argument$.$ Thus
\vskip-2pt\noi
$$
\|\vphi\otimes\eta\|_{*[\alpha,\alpha]}\le\|\vphi\|\,\|\eta\|.
$$
\vskip4pt\noi
An arbitrary element $\Bbbk$ in
${\B[\X,\V]^*\otimes\B[\Y,\W]^*}\sse
{\Le\big[\B[\X,\V]\otimes_{[\alpha,\alpha]}\B[\Y,\W],\FF\kern.5pt\big]}$
is represented by a finite sum $\Bbbk={\sum_k\vphi_k\otimes\eta_k}$ of single
tensors with ${\vphi_k\in\B[\X,\V]^*}$ and ${\eta_k\in\B[\Y,\W]^*}.$ Take
the induced uniform norm on
${\B\big[\B[\X,\V]\otimes_{[\alpha,\alpha]}\B[\Y,\W],\FF\kern.5pt\big]}$, say,
${\|\cdot\|_{*[\alpha,\alpha]}}.$ The above displayed inequality ensures that
\begin{eqnarray*}
\|\Bbbk\|_{*[\alpha,\alpha]}
&\kern-6pt=\kern-6pt&
\sup_{\|\LL\|_{[\alpha,\alpha]}\le1}
\Big|\Big({\sum}_k\vphi_k\otimes\eta_k\Big)\LL\Big|                       \\
&\kern-6pt\le\kern-6pt&
\sup_{\|\LL\|_{[\alpha,\alpha]}\le1}
{\sum}_k\|\vphi_k\otimes\eta_k\|_{*[\alpha,\alpha]}\|\LL\|_{[\alpha,\alpha]}
\le{\sum}_k\|\vphi_k\|\,\|\eta_k\|,
\end{eqnarray*}
which is finite as the sum is finite$.$ Thus
${\Bbbk\in\kern-1pt\B\big[\B[\X,\V]\otimes_{[\alpha,\alpha]}
\kern-1pt\B[\Y,\W],\FF\kern.5pt]\kern.5pt\big]}$,
\hbox{and hence}
$$
\B[\X,\V]^*\otimes\B[\Y,\W]^*
\sse\B\big[\B[\X,\V]\otimes_{[\alpha,\alpha]}\B[\Y,\W],\FF\kern.5pt\big]
=(\B[\X,\V]\otimes_{[\alpha,\alpha]}\B[\Y,\W]\kern.5pt)^*.
$$
So equip ${\B[\X,\V]^*\otimes\B[\Y,\W]^*}$ with the norm
${\|\cdot\|_{*[\alpha,\alpha]}}$ on
${([\B[\X,\V]\otimes_{[\alpha,\alpha]}\B[\Y,\W]\kern.5pt)^*}$ and set
${\B[\X,\V]^*\otimes_{*[\alpha,\alpha]}\B[\Y,\W]^*
=(\B[\X,\V]^*\otimes\B[\Y,\W]^*,\|\cdot\|_{*[\alpha,\alpha]})}$.        \qed

\section{Main Result}

Theorem 4.1 shows that ${\|\cdot\|_{[\alpha,\alpha]}}$ is a reasonable
crossnorm on ${\B[\X,\V]\otimes\B[\Y,\W]}$, when inherited from
${\B[\X\otimes_\alpha\!\Y,\,\V\otimes_\alpha\!\W]}$, thus extending
Proposition 2.2 from continuous linear functionals to arbitrary
continuous linear transformations$.$

\vskip6pt
The proof of Theorem 4.1 is especially tailored to prompt the question that
closes the paper, and also to support the statement of Corollary 4.4.

\vskip3pt
\theorem{4.1}
{\it If\/ ${\|\cdot\|_\alpha}$ is a uniform crossnorm, then
$$
\|\LL\|_\ve\le\|\LL\|_{[\alpha,\alpha]}\le\|\LL\|_\we         \leqno{\rm(a)}
$$
for every\/ ${\LL\in\B[\X,\V]\otimes_{[\alpha,\alpha]}\B[\Y,\W]}$,
where\/ ${\X,\Y,\V,\W}$ are arbitrary normed spaces and\/
${\|\cdot\|_{[\alpha,\alpha]}}$ is the associated induced uniform norm on\/
${\B[\X\otimes_\alpha\!\Y,\V\otimes_\alpha\!\W]}.$
Hence
\vskip6pt\noi
{\rm(b)}
$\;\;{\|\cdot\|_{[\alpha,\alpha]}}$ is a reasonable crossnorm on}\/
${\B[\X,\V]\otimes\B[\Y,\W]}$.

\proof
(a)
Take an arbitrary element
${\LL\in\B[\X,\V]\otimes\B[\Y,\W]}
\sse{\B[\X\otimes_\alpha\!\Y,\V\otimes_\alpha\!\W]}$
(according to Proposition 3.1) and let ${\sum_jA_j\otimes B_j}$ be any
finite-sum representation of $\LL$ in terms of single tensors in
${\B[\X,\V]\otimes\B[\Y,\W]}.$ Regard $\LL$ as a transformation in
${\B[\X\otimes_\alpha\!\Y,\V\otimes_\alpha\!\W]}.$ Set
${\|\LL(\digamma)\|_{\V\otimes_\alpha\W}=\|\LL(\digamma)\|_\alpha}$ and
${\|\digamma\|_{\X\otimes_\alpha\Y}=\|\digamma\|_\alpha}$ for
$\digamma$ in ${\X\!\otimes_\alpha\!\Y}.$ With ${\|\cdot\|_{[\alpha,\alpha]}}$
standing for the induced uniform norm on
${\B[\X\!\otimes_\alpha\!\Y,\V\kern-1pt\otimes_\alpha\!\W]}$,
\begin{eqnarray*}
\|\LL(\digamma)\|_\alpha
&\kern-6pt=\kern-6pt&
\Big\|\Big({\sum}_jA_j\otimes B_j\Big)\digamma\Big\|_\alpha
\le{\sum}_j\|(A_j\otimes B_j)\digamma\|_\alpha                            \\
&\kern-6pt\le\kern-6pt&
{\sum}_j\|A_j\otimes B_j\|_{[\alpha,\alpha]}
\|\digamma\|_\alpha
=\Big({\sum}_j\|A_j\|\,\|B_j\|\Big)\|\digamma\|_\alpha
\end{eqnarray*}
for every ${\digamma\!\in\kern-1pt\X\otimes\Y}$. As
${\|\LL\|_{[\alpha,\alpha]}}={\sup_{\|\digamma\|_\alpha\le1}
\|\LL(\digamma)\|_\alpha}$ we get
$$
\|\LL\|_{[\alpha,\alpha]}\le{\sum}_j\|A_j\|\,\|B_j\|.
$$
Since the above inequality holds for every representation
${\sum_jA_j\otimes B_j}$ of $\LL$,
$$
\|\LL\|_{[\alpha,\alpha]}
\le\inf_{\LL=\Sigma_jA_j\otimes B_j}{\sum}_j\|A_j\|\,\|B_j\|=\|\LL\|_\we,
$$
where ${\|\cdot\|_\we}\kern-1pt$ is the projective norm on
${\B[\X\kern-1pt,\V]\otimes\B[\Y,\W]}.$ On the other hand, take an arbitrary
$\Bbbk={\sum_k\vphi_k\otimes\eta_k}$ in
${\B[\X,\V]^*\otimes\B[\Y,\W]^*
\sse(\B[\X,\V]\otimes_{[\alpha,\alpha]}\B[\Y,\W])^*}$
(according to Proposition 3.2)$.$ Regard $\Bbbk$ as a functional in
$(\B[\X,\V]\otimes_{[\alpha,\alpha]}\B[\Y,\W])^*.$ With
${\|\cdot\|_{*[\alpha,\alpha]}}$ being the induced uniform norm in
${\B\big[\B[\X,\V]\otimes_{[\alpha,\alpha]}\B[\Y,\W],\FF\kern.5pt\big]}$,
$$
\|\Bbbk\|_{*[\alpha,\alpha]}
=\sup_{\|\LL\|_{[\alpha,\alpha]}\le1}|\Bbbk(\LL)|.
\;\quad\hbox{Dually,}\quad\;
\|\LL\|_{[\alpha,\alpha]}=\sup_{\|\Bbbk\|_{*[\alpha,\alpha]}\le1}|\Bbbk(\LL)|.
$$
Since for every
${\Bbbk=\sum_k\vphi_k\otimes\eta_k\in\B[\X,\V]^*\otimes\B[\Y,\W]^*}$
the value of ${\Bbbk(\LL)\in\FF}$ is
$$
\Bbbk(\LL)
=\Big({\sum}_k\vphi_k\otimes\eta_k\Big){\sum}_jA_j\otimes B_j
={\sum}_{k,j}\vphi_k(A_j)\otimes\eta_k(B_j)
={\sum}_{k,j}\vphi_k(A_j)\,\eta_k(B_j)
$$
for every ${\LL=\sum_jA_j\otimes B_j\in\B[\X,\V]\otimes\B[\Y,\W]}\kern.5pt$,
we get
\begin{eqnarray*}
\|\LL\|_{[\alpha,\alpha]}
&\kern-6pt=\kern-6pt&
\sup_{\|\Bbbk\|_{*[\alpha,\alpha]}\le1}|\Bbbk(\LL)|
=\sup_{\|\Sigma_k\vphi_k\otimes\eta_k\|_{*[\alpha,\alpha]}\le1}
\Big|{\sum}_{k,j}\vphi_k(A_j)\,\eta_k(B_j)\Big|                           \\
&\kern-6pt\ge\kern-6pt&
\sup_{\|\vphi\otimes\eta\|_{*[\alpha,\alpha]}\le1}
\Big|{\sum}_j\vphi(A_j)\,\eta(B_j)\Big|
\ge\sup_{\|\vphi\|\kern1pt\|\eta\|\le1}
\Big|{\sum}_j\vphi(A_j)\,\eta(B_j)\Big|,
\end{eqnarray*}
because ${\|\vphi\otimes\eta\|_{*[\alpha,\alpha]}}\le{\|\vphi\|\,\|\eta\|}$
for every ${\vphi\in\B[\X,\V]^*}$ and ${\eta\in\B[\Y,\W]^*}$ according to
Proposition 3.2$.$ Hence
$$
\|\LL\|_{[\alpha,\alpha]}
\ge\sup_{\|\vphi\|\kern1pt\|\eta\|\le1}\Big|{\sum}_j\vphi(A_j)\,\eta(B_j)\Big|
\ge\sup_{\|\vphi\|\le1,\,\|\eta\|\le1}\Big|{\sum}_j\vphi(A_j)\,\eta(B_j)\Big|
=\|\LL\|_\ve,
$$
where ${\|\cdot\|_\ve}$ is the injective norm on
${\B[\X,\V]\otimes\B[\Y,\W]}.$ Consequently,
$$
\|\LL\|_\ve\le\|\LL\|_{[\alpha,\alpha]}\le\|\LL\|_\we
$$
for every ${\LL\in\B[\X,\V]\otimes_{[\alpha,\alpha]}\B[\Y,\W]}$.

\vskip6pt\noi
(b)
Thus, according to Proposition 2.1, the induced uniform norm
${\|\cdot\|_{[\alpha,\alpha]}}$ becomes a {\it reasonable crossnorm}\/ on the
tensor product space ${\B[\X,\V]\otimes\B[\Y,\W]}$.                      \qed

\vskip9pt
{\it Particular case}\/$.$
Set ${\V\kern-1pt=\X}$, ${\W\kern-1pt=\!\Y}$ and write $\B[\X]$, $\B[\Y]$,
and ${\B[\X\otimes_\alpha\!\Y]}$ for ${\B[\X,\X]}$, $\,{\B[\Y,\Y]}$, and
$\,{\B[\X\otimes_\alpha\!\Y,\X\otimes_\alpha\!\Y]}$, respectively$.$ Thus
Theorem 4.1 yields the result stated in the Abstract on tensor products of
algebra of operators.

\vskip6pt
It is worth noticing that a first application of the inequalities in part (a)
of Theorem 4.1 is the assignment of a reasonable crossnorm to the tensor
product space ${\B[\X,\V]\otimes\B[\Y,\W]}$ associated with a single uniform
crossnorm as in part (b).

\vskip3pt
\remark{4.2}
If ${\|\cdot\|_\alpha}$ is a uniform crossnorm, then (by definition)
$$
\B[\X,\V]\otimes\B[\Y,\W]\sse\B[\X\otimes_\alpha\!\Y,\V\otimes_\alpha\!\W].
$$
\vskip-4pt\noi
As before, let
$$
\hbox{$\|\cdot\|_{[\alpha,\alpha]}$ be the induced uniform norm on
$\B[\X\otimes_\alpha\!\Y,\V\otimes_\alpha\!\W]$}.
$$
\vskip2pt\noi
Take an arbitrary
$\LL=\kern-1pt\sum_jA_j\otimes B_j\kern-1pt\in\B[\X,\V]\otimes\B[\Y,\W]
\sse\Le[\X\otimes\Y,\V\otimes\W].$
Since ${\|\cdot\|_\ve\le\|\cdot\|_\alpha\le\|\cdot\|_\we}$
(cf$.$ Proposition 2.1),
$$
{\sup}_{\digamma\ne0}
\smallfrac{\|\LL(\digamma)\|_\ve}{\|\digamma\|_\we}
\le{\sup}_{\digamma\ne0}
\smallfrac{\|\LL(\digamma)\|_\alpha}{\|\digamma\|_\alpha}
\le{\sup}_{\digamma\ne0}
\smallfrac{\|\LL(\digamma)\|_\we}{\|\digamma\|_\ve}.              \eqno{(*)}
$$
Thus set
${\|\LL\|_{[\we,\ve]}}
={\sup_{\digamma\ne0}\smallfrac{\|\LL(\digamma)\|_\ve}{\|\digamma\|_\we}}$
which is finite as
${\|\LL\|_{[\alpha,\alpha]}}
={\sup_{\digamma\ne0}
\smallfrac{\|\LL(\digamma)\|_\alpha}{\|\digamma\|_\alpha}}$
is a norm, which is enough to ensure that 
$$
\hbox{${\|\cdot\|_{[\we,\ve]}}$
is the induced uniform a norm on ${\B[\X\otimesw\Y,\V\otimesv\W]}$},
$$
and that ${\LL\in\B[\X\otimesw\Y,\V\otimesv\W]}.$ (Note:
${\sup_{\digamma\ne0}\smallfrac{\|\LL(\digamma)\|_\we}{\|\digamma\|_\ve}}$ may
not be finite$.$) So
$$
\B[\X\otimes_\alpha\!\Y,\V\otimes_\alpha\!\W]
\sse\B[\X\otimesw\Y,\V\otimesv\W]
\quad\hbox{where}\quad
\|\cdot\|_{[\we,\ve]}\le\|\cdot\|_{[\alpha,\alpha]}
$$
by $(*)$, and since ${\LL\in\B[\X\otimesw\Y,\V\otimesv\W]}$ we also get
$$
\B[\X,\V]\otimes\B[\Y,\W]\sse\B[\X\otimesw\Y,\V\otimesv\W].
$$
Therefore we may equip ${\B[\X,\V]\otimes\B[\Y,\W]}$ with both induced
uniform norms, namely, ${\|\cdot\|_{[\ve,\we]}}$ and
${\|\cdot\|_{[\alpha,\alpha]}}$, and so we may consider the normed spaces
$$
\B[\X,\V]\otimes_{[\alpha,\alpha]}\B[\Y,\W]
\sse\B[\X,\V]\otimes_{[\we,\ve]}\B[\Y,\W]
\quad\hbox{with}\quad
\|\cdot\|_{[\we,\ve]}\le\|\cdot\|_{[\alpha,\alpha]}.
$$
\vskip-2pt

\vskip3pt
\remark{4.3}
For the induced uniform norm ${\|\cdot\|_{[\we,\ve]}}$ on
${\B[\X\otimesw\Y,\V\otimesv\W]}$,
$$
\|\LL\|_\ve\le\|\LL\|_{[\we,\ve]}\le\|\LL\|_\we
$$
for every ${\LL\in\B[\X,\V]\otimes_{[\we,\ve]}\B[\Y,\W]}$, and hence
$$
\hbox{${\|\cdot\|_{[\we,\ve]}}$ is a reasonable crossnorm on
${\B[\X,\V]\otimes\B[\Y,\W]}$}.
$$
Indeed, it can be verified that ${\|\cdot\|_\ve\le\|\cdot\|_{[\we,\ve]}}.$
Also, ${\|\cdot\|_{[\we,\ve]}\le\|\cdot\|_{[\alpha,\alpha]}}$ for an arbitrary
uniform crossnorm ${\|\cdot\|_\alpha}$ according to Remark 4.2$.$ Moreover,
Theorem 4.1(a) says that ${\|\cdot\|_{[\alpha,\alpha]}\le\|\cdot\|_\we}.$
Summing up,
$$
\|\cdot\|_\ve\le\|\cdot\|_{[\we,\ve]}
\le\|\cdot\|_{[\alpha,\alpha]}\le\|\cdot\|_\we.
$$
\vskip0pt\noi
Again, Proposition 2.1 ensures that the induced uniform norm
${\|\cdot\|_{[\we,\ve]}}$ is a reasonable crossnorm on the
${\B[\X,\V]\otimes\B[\Y,\W]}$.

\vskip9pt
The injective and projective norms ${\|\cdot\|_\ve}$ and ${\|\cdot\|_\we}$,
being uniform crossnorms, act on every tensor product space and are the least
and the greatest uniform crossnorm on every tensor product space (cf$.$
Proposition 2.1)$.$ In particular, on
${\B[\X,\V]\kern-1pt\otimes\kern-1pt\B[\Y,\W]}$ (as in the above displayed
inequalities)$.$ Also, ${\|\cdot\|_{[\ve,\ve]}}\kern-1pt$ and
${\|\cdot\|_{[\we,\we]}}$ are reasonable crossnorms on
${\B[\X,\V]\otimes\B[\Y,\W]}$ by Theorem 4.1 because ${\|\cdot\|_\ve}$ and
${\|\cdot\|_\we}$ are uniform crossnorm$.$ Therefore
$$
\|\LL\|_\ve\le\|\LL\|_{[\ve,\ve]}
\quad\hbox{and}\quad
\|\LL\|_{[\we,\we]}\le\|\LL\|_\we.
$$
\vskip-2pt

\vskip6pt\noi
{\it Question}\/$.$
Does ${\|\LL\|_{[\alpha,\alpha]}}$ lie in between, so that
${\|\LL\|_{[\ve,\ve]}}$ is the least and ${\|\LL\|_{[\we,\we]}}$ is the
greatest reasonable crossnorm on ${\B[\X,\V]\otimes\B[\Y,\W]}$ that are
inherited from ${\B[\X\otimes_\alpha\!\Y,\V\otimes_\alpha\!\W]}$ for arbitrary
uniform crossnorms ${\|\cdot\|_\alpha}\kern.5pt$? (Note that there are
reasonable crossnorms on the tensor product space ${\B[\X,\V]\otimes\B[\Y,\W]}$
that are not restrictions of the induced uniform norm on
${\B[\X\otimes_\alpha\!\Y,\V\otimes_\alpha\!\W]}$ for any uniform crossnorm
${\|\cdot\|_\alpha}$, as is the case of ${\|\LL\|_\ve}$ and ${\|\LL\|_\we}$.)

\vskip6pt
Corollary 4.4 shows how ${\|\LL\|_{[\alpha,\alpha]}}$ naturally fits between
${\|\LL\|_{[\ve,\ve]}}$ and ${\|\LL\|_{[\we,\we]}}$, giving a first estimate
to the above a question$.$ (Compare the result in Corollary 4.4 below with
$(*)$ in Remark 4.2.) Consider the setup in the proof of \hbox{Theorem 4.1.}

\vskip3pt
\corollary{4.4}
{\it If\/ ${\|\cdot\|_\alpha}$ is a uniform crossnorm, then}\/
$$
\inf_{0\ne\digamma\in\X\otimes\Y}
\smallfrac{\|\digamma\|_\ve}{\|\digamma\|_\alpha}\;
\|\LL\|_{[\ve,\ve]}
\le\|\LL\|_{[\alpha,\alpha]}
\le\sup_{0\ne\digamma\in\X\otimes\Y}
\smallfrac{\|\digamma\|_\we}{\|\digamma\|_\alpha}\;
\|\LL\|_{[\we,\we]}.
$$

\proof
Consider the uniform crossnorms norms ${\|\cdot\|_\ve}$, ${\|\cdot\|_\alpha}$,
and ${\|\cdot\|_\we}$ acting on ${\X\otimes\Y}$ and on ${\V\otimes\W}.$
According to Proposition 2.1
${\|\cdot\|_\ve}\!\le{\|\cdot\|_\alpha}\!\le{\|\cdot\|_\we}.$ Regarding the
setup in the proof of Theorem 4.1, take ${\LL(\digamma)\in\V\otimes\W}$ for an
arbitrary ${\digamma\!\in\X\otimes\Y}.$ First consider the inequality
${\|\LL(\digamma)\|_\alpha\le\|\LL(\digamma)\|_\we}$ so that
$$
\|\LL\|_{[\alpha,\alpha]}
=\sup_{\digamma\ne0}\smallfrac{\|\LL(\digamma)\|_\alpha}{\|\digamma\|_\alpha}
\le\sup_{\digamma\ne0}\smallfrac{\|\LL(\digamma)\|_\we}{\|\digamma\|_\alpha}
\le\|\LL\|_{[\we,\we]}\,
\sup_{\digamma\ne0}\smallfrac{\|\digamma\|_\we}{\|\digamma\|_\alpha}.
$$
Next take the inequality ${\|\LL(\digamma)\|_\ve\le\|\LL(\digamma)\|_\alpha}.$
Similarly,
$$
\|\LL\|_{[\ve,\ve]}
\le\|\LL\|_{[\alpha,\alpha]}\,
\sup_{\digamma\ne0}\smallfrac{\|\digamma\|_\alpha}{\|\digamma\|_\ve};
\quad\hbox{equivalently,}\quad
\inf_{\digamma\ne0}\smallfrac{\|\digamma\|_\ve}{\|\digamma\|_\alpha}\,
\|\LL\|_{[\ve,\ve]}
\le\|\LL\|_{[\alpha,\alpha]},
$$
with
${\inf_{\digamma\ne0}\smallfrac{\|\digamma\|_\ve}{\|\digamma\|_\alpha}}
={\big(\sup_{\digamma\ne0}
\smallfrac{\|\digamma\|_\alpha}{\|\digamma\|_\ve}\big)^{-1}}.$    \qed

\vskip3pt
\remark{4.5}
It is not usual to consider more than one uniform crossnorm, one equipping the
domain and the other equipping the codomain, of a tensor product ${A\otimes B}$
of bounded linear transformations --- see, e.g., \cite[Section 6.1]{Rya} and
\cite[Section 1.2]{DFS}$.$ In fact, a uniform crossnorm ${\|\cdot\|_\alpha}$
is a reasonable crossnorm which is supposed to be assigned to every tensor
product for arbitrary normed spaces making ${A\otimes B}$ continuous when
acting from ${\X\otimes_\alpha\!\Y}$ to ${\V\otimes_\alpha\!\W}$; both tensor
product spaces equipped with the same norm (cf$.$ Section 3)$.$ Therefore, if
one takes another reasonable crossnorm ${\|\cdot\|_\beta}$, also supposed to
be assigned to \hbox{every} tensor product space of arbitrary normed spaces,
and requires that ${A\in\B[\X,\V]}$ and ${B\in\B[\Y,\W]}$ imply
${A\otimes B\in\B[\X\otimes_\alpha\!\Y,\V\otimes_\beta\!\W]}$, then it is
readily verified that
${\|A\|\kern1pt\|B\|}\le{\|A\otimes B\|}={\|A\otimes B\|_{[\alpha,\beta]}}.$
A new concept of a {\it jointly uniformly crossnorms}\/ would require a new
definition to make the above inequality an identity$.$ Also, under a new
definition of a jointly uniform crossnorm, Propositions 3.1 and 3.2 would
require a new restatement$.$ (Extending Proposition 3.1 to such a new setup
seems to be a simple task, but Proposition 3.2 would perhaps require some
additional, possibly nontrivial, arguments)$.$ Although we will not proceed
along this line --- it goes beyond the purpose of the present paper --- it
seems that a possible version of the first part of Theorem 4.1 involving a
pair of distinct {\it jointly uniform crossnorm}\/ might lead to a promising
further research.

\section{Final Remark}

Suppose a pair of normed spaces ${(\X,\Y)}$ is such that
${\sup_{\,0\ne\digamma\in\X\otimes\Y}\frac{\|\digamma\|_\we}{\|\digamma\|_\ve}
\!<\!\infty}$
when their tensor product space ${\X\otimes\Y}$ is equipped with the injective
norm ${\|\cdot\|_\ve}$ and with the projective norm ${\|\cdot\|_\we}.$ (Trivial
example: if $\X$ and $\Y$ are finite-dimen\-sional, where all norms are
equivalent, so that ${\X\otimesv\Y\cong\X\otimesw\Y}$ --- here $\cong$ means
topological isomorphism)$.$ In such a case (when the above supremum is finite),
set
${\sup_{\,0\ne\digamma\in\X\otimes\Y}\frac{\|\digamma\|_\we}{\|\digamma\|_\ve}
=\gamma}$,
and the injective and projective norms become {\it equivalent uniform
crossnorms}\/ on ${\X\otimes\Y}$ and so is any uniform crossnorm
\hbox{${\|\cdot\|_\alpha}$ since}
$$
{\|\cdot\|_\ve\le\|\cdot\|_\alpha\le\|\cdot\|_\we\le\gamma\|\cdot\|_\ve},
$$
and so
$\sup_{\digamma\ne0}\frac{\|\digamma\|_\we}{\|\digamma\|_\alpha}\!<\!\infty$
and
$\sup_{\digamma\ne0}\frac{\|\digamma\|_\alpha}{\|\digamma\|_\ve}\!<\!\infty.$
It is attributed to Grothendieck the origin of the question whether
${\X\hotimesv\Y}\cong{\X\hotimesw\Y}$ holds for some pair of
infinite-dimensional Banach spaces $\X$ and $\Y$ (see \cite[p.181]{Pis0})$.$
A solution was given by Pisier in \cite[Theorem 3.2(b)]{Pis0} where it was
exhibited a separable infinite-dimensional \hbox{Banach} space $\Pe\!$, now
called {\it Pisier space}\/, such that
${\Pe\hotimesv\Pe\cong\Pe\hotimesw\kern-1pt\Pe}$ (here $\cong$ means isometric
isomorphism), which shows in addition that all reasonable crossnorms (and all
uniform crossnorm) on ${\Pe\hotimes\Pe}$ (and so on ${\Pe\otimes\Pe}$) are
isomorphically equivalent, where in this case ${\gamma=1}$ and, consequently,
for ${\LL\in\B[\Pe]\otimes\B[\Pe]}$,
$$
\|\LL\|_{[\ve,\ve]}=\|\LL\|_{[\alpha,\alpha]}=\|\LL\|_{[\we,\we]},
$$
with respect to the setup in the proof of Theorem 4.1$.$ Corollary 4.4 gives
just a first estimate to the question whether the equal signs ``$=$'' in the
above equation can be replaced by ``$\le$'' (perhaps weighted with positive
constants) for an arbitrary transformation ${\LL\in\B[\X,\V]\otimes\B[\Y,\W]}$
acting on arbitrary (or on specific classes of) normed spaces ${\X,\Y,\V,\W}$.

\section*{Acknowledgment}

The author thanks a referee for suggestions to improve the paper which were
included in Remarks 4.2 and 4.3; especially for calling his attention to the
inequality that opens the argument of Remark 4.3.

\bibliographystyle{amsplain}

\begin{thebibliography}{10}

\bibitem{Arv}
W. Arveson,
{\it What is an operator space\/?}\/,
(2008), 1--7, available at
\vskip0pt\noi
https://math.berkeley.edu/$\sim$arveson/Dvi/opSpace.pdf

\bibitem{BP}
D.P. Blecher and V.I. Paulsen,
{\it Tensor products of operator spaces}\/,
J. Funct. Anal. {\bf 99} (1991), 262--292.

\bibitem{DF}
A. Defant and K. Floret,
{\it Tensor Norms and Operator Ideals}\/,
North-Holland, Amsterdam, 1993.

\bibitem{DFS}
J. Diestel, J.H. Fourie, and J. Swart,
{\it The Metric Theory of Tensor Products -- Grothendieck's R\'esum\'e
Revisited}\/,
American Mathematical Society, Providence, 2008.

\bibitem{Jar}
H. Jarchow,
{\it Locally Convex Spaces}\/,
B.G. Teubner, Stuttgart, 1981.

\bibitem{EOT}
C.S. Kubrusly,
{\it The Elements of Operator Theory}\/,
2nd edn. Birkh\"auser-Springer, New York, 2011.

\bibitem{Kub}
C.S. Kubrusly,
{\it Algebraic tensor products revisited$:$ Axiomatic approach}\/,
Bull. Malays. Math. Sci. Soc. {\bf 44} (2021), 2335--2355.

\bibitem{Pis0}
G. Pisier,
{\it Counterexamples to a conjecture of Grothendieck}\/,
Acta Math. {\bf 151} (1983), 181--208.

\bibitem{Pis1}
G. Pisier,
{\it Introduction to Operator Space Theory}\/,
Cambridge University Press, Cambridge, 2003.

\bibitem{Pis2}
G. Pisier,
{\it Tensor Products of C*-Algebras and Operator Spaces}\/,
Cambridge University Press, Cambridge, 2020.

\bibitem{Rya}
R.A. Ryan,
{\it Introduction to Tensor Products of Banach Spaces}\/,
Springer, London, 2002.

\bibitem{Sim}
B. Simon,
{\it Uniform crossnorms}\/,
Pacific J. Math. {\bf 46} (1973), 555--560.

\end{thebibliography}

\end{document}